\input amstex
\documentstyle{amsppt}
\pagewidth{6.2in}\vsize8.5in\parindent=6mm
\parskip=3pt\baselineskip=14pt\tolerance=10000\hbadness=500
\document
\topmatter
\title Note on the Prime Number Theorem
\endtitle
\author Yong-Cheol Kim \endauthor
\abstract We survey the classical results on the prime number
theorem.
\endabstract
\address Department of Mathematics, Korea University, Seoul 136-701, Korea
\endaddress
\endtopmatter

\define\inn#1#2{\langle#1,#2\rangle}

\define\lcontr{\rfloor}
\define\lco#1#2{{#1}\lcontr{#2}}
\define\lcoi#1#2{\imath({#1}){#2}}
\define\rco#1#2{{#1}\rcontr{#2}}

\redefine\Rea{{\text{\rm Re}}} \redefine\arg{{\text{\rm arg }}}
\redefine\Ima{{\text{\rm Im}}}

\define\ap{\alpha}             
\define\bt{\beta}
\define\gm{\gamma}             
\define\dt{\delta}             
\define\vep{\varepsilon}
\define\zt{\zeta}
\define\th{\theta}             \define\Th{\Theta}

\define\ld{\lambda}            \define\Ld{\Lambda}
\define\sm{\sigma}             

\define\om{\omega}             \define\Om{\Omega}
            \define\iy{\infty}
\define\lt{\left}            \define\rt{\right}
\define\f{\frac}

\define\fU{{\frak U}}


\define\fm{{\frak m}}
\define\fn{{\frak n}}
\define\fo{{\frak o}}
\define\fp{{\frak p}}

\define\fu{{\frak u}}

\define\BB{{\Bbb B}}
\define\BC{{\Bbb C}}

\define\BN{{\Bbb N}}

\define\BR{{\Bbb R}}

\define\BZ{{\Bbb Z}}

\define\cO{{\Cal O}}




\define\s{\setminus}         
            
\define\pa{\partial}        \define\fd{\fallingdotseq}
       
     \define\ds{\dsize}
\redefine\li{{\text{\rm li }}}

In this chapter, we are very interested in the asymptotic behavior
of a single number theoretic function $\pi(n)$ which counts all
prime numbers between $1$ and $n$, or $\pi(x)$ which is extended
to $\BR$ and defined by $$\pi(x)=\sum_{p\le x}\,1.$$ It is
well-known that Euclid showed that $$\lim_{x\to\iy}\pi(x)=\iy\,;$$
that is, there exist infinitely many prime numbers.

\proclaim{Proposition 5.1} There exists a constant $c>0$ such that
$$\pi(x)\ge c\cdot\ln\ln x.$$ \endproclaim

\noindent{\it Proof.} First of all, we prove that if $p_n$ is the
$n$th prime number then we have that $$p_n\le 2^{2^{n-1}}.$$ Since
there must be some $p_{n+1}$ dividing the number $p_1 p_2\cdots
p_n -1$ and not exceeding it, it follows from the induction step
that
$$p_{n+1}\le 2^{2^0}\,2^{2^1}\cdots
2^{2^{n-1}}=2^{2^0+2^1+\cdots+2^{n-1}}\le 2^{2^n}.$$ If $x\ge 2$
is some real number, then we select the largest natural number $n$
satisfying $\,2^{2^{n-1}}\le x,$ so that we have that
$\,2^{2^n}>x$. Hence we conclude that
$$\pi(x)\ge n\ge\f{1}{\ln 2}\cdot\ln\lt(\f{\ln x}{\ln 2}\rt)\ge\f{1}{\ln 2}\cdot\ln\ln x.\,\,\,\,\qed$$

\proclaim{Proposition 5.2} There exists a constant $c>0$ such that
$$\pi(x)\ge c\cdot\ln x$$ for all sufficiently large $x$.\endproclaim

\noindent{\it Proof.} Since each square-free integer $n\le x$ can
be only be divided by $p_1,p_2,\cdots,p_{\pi(x)}$, $n$ can be
written uniquely as $$n=\prod_{k=1}^{\pi(x)} p_k^{\ap_k}$$ where
$\ap_k$ takes only the values $0$ or $1$. Thus there are at most
$2^{\pi(x)}$ square-free integers $n\le x$. From Corollary 4.2.21,
we see that the density of the square-free integers tends to
$6/\pi^2\,$; that is, the number of square-free numbers $n\le x$
grows asymptotically to $6x/\pi^2$. This implies that there is
some constant $c_0<6/\pi^2$ such that $$c_0\cdot x\le 2^{\pi(x)}$$
for all sufficiently large $x$. Hence we complete the proof. \qed

Neither of Proposition 5.1 and Proposition 5.2 describes the
asymptotic behavior of $\pi(x)$ quite well. Long time ago,
Legendre and Gauss conjectured that $$\pi(x)\sim\f{x}{\ln x}.$$
The truth of this assertion is the core of the prime number
theorem. For more delicate description of $\pi(x)$, we consder the
integral logarithm function $\li x$ defined as the Cauchy
principal value integral $$\li x=\int_0^x\f{1}{\ln
t}\,dt=\lim_{\vep\to\iy}\lt(\int_0^{1-\vep}\f{1}{\ln
t}\,dt+\int_{1+\vep}^x\f{1}{\ln t}\,dt\rt).$$ It follows from de
l'Hospital's rule that $$\lim_{x\to\iy}\f{\li x}{\ds\f{x}{\ln
x}}=\lim_{x\to\iy}\f{\ds\f{1}{\ln x}}{\ds\f{1}{\ln x}-\f{1}{\ln^2
x}}=1.$$ Thus we obtain the asymptotic behavior of $\li x$ as
follows; $$\li x\sim\f{x}{\ln x}.$$ Hence the asymptotic relation
$\,\pi(x)\sim\li x$ is called the prime number theorem. In fact,
Gauss conjectured that $\li x$ describes $\pi(x)$ even better than
$x/\ln x$.

\proclaim{Lemma 5.3} $(a)$ $\ds\sum_{n\le
x}\Ld(n)\lt[\f{x}{n}\rt]=x\ln x-x+\cO(\ln x).$

$(b)$ $\ds\sum_{n\le
x}\Ld(n)\lt(\lt[\f{x}{n}\rt]-2\lt[\f{x}{2n}\rt]\rt)=x\ln 2+\cO(\ln
x).$ \endproclaim

\noindent{\it Proof.} (a) By the definition of the Mangoldt
function, we have that $$\sum_{n\le x}\ln n=\sum_{n\le
x}\sum_{m|n}\Ld(m)=\sum_{m\le x}\Ld(m)\sum_{n\le x : m|n}
1=\sum_{m\le x}\Ld(m)\lt[\f{x}{m}\rt].$$ Thus it follows from
Proposition 4.2.3[the Euler's sum formula] that $$ \sum_{n\le
x}\Ld(n)\lt[\f{x}{n}\rt]=\sum_{n\le x}\ln n =\int_1^x\ln
t\,dt+\cO(\ln x) =x\ln x-x+\cO(\ln x). $$

(b) By applying (a) and the fact that $\ds\sum_{x/2<n\le
x}\Ld(n)\lt[\f{x}{2n}\rt]=0,$ we obtain that $$ \split \sum_{n\le
x}\Ld(n)\lt(\lt[\f{x}{n}\rt]-2\lt[\f{x}{2n}\rt]\rt)&=\sum_{n\le
x}\Ld(n)\lt[\f{x}{n}\rt]-2\sum_{n\le
x/2}\Ld(n)\lt[\f{x}{2n}\rt]-2\sum_{x/2<n\le
x}\Ld(n)\lt[\f{x}{n}\rt] \\ &=x\ln
x-x-2\lt(\f{x}{2}\,\ln\f{x}{2}-\f{x}{2}\rt)+\cO(\ln x) \\
&=x\ln 2+\cO(\ln x).\endsplit $$ Hence we complete the proof.\qed

\proclaim{Theorem 5.4[Chebyshev's Theorem]} There exist two
constants $c_1>0$ and $c_2>0$ such that
$$c_1\cdot\f{x}{\ln x}\le\pi(x)\le c_2\cdot\f{x}{\ln x}$$ for all
sufficiently large $x$. \endproclaim

\noindent{\it Proof.} Since $\ds [\ap]-2\lt[\f{\ap}{2}\rt]$ is
always an integer and satisfies the following inequality
$$-1=\ap-1-2\,\f{\ap}{2}<[\ap]-2\lt[\f{\ap}{2}\rt]<\ap-2\lt(\f{\ap}{2}-1\rt)=2,$$
we see that $$0\le [\ap]-2\lt[\f{\ap}{2}\rt]\le 1.\tag{5.1}$$ Thus
by (5.1) and (b) of Lemma 5.3 we have that $$ \split x\ln
2+\cO(\ln x)&=\sum_{n\le
x}\Ld(n)\lt(\lt[\f{x}{n}\rt]-2\lt[\f{x}{2n}\rt]\rt) \\
&\le\sum_{n\le x}\Ld(n)=\sum_{p\le x}\lt[\f{\ln x}{\ln p}\rt]\,\ln
p \\ &\le\ln x\,\sum_{p\le x}\,1=\pi(x)\,\ln x, \endsplit $$ and
so we can get the first inequality by dividing by $\ln x$. For the
second inequality, we observe that $$ \split \pi(x)\ln
x-\pi\lt(\f{x}{2}\rt)\ln\f{x}{2}&=\ln\f{x}{2}\,\lt(\pi(x)-\pi\lt(\f{x}{2}\rt)\rt)+\pi(x)\,\ln2
\\ &=\ln\f{x}{2}\,\lt(\pi(x)-\pi\lt(\f{x}{2}\rt)\rt)+\cO(x) \\
&=\cO\lt(\sum_{x/2<p\le x}\ln p +x\rt) \\
&=\cO\lt(\sum_{x/2<n\le x}\Ld(n)\cdot (1-0)+x\rt) \\
&=\cO\lt(\sum_{x/2<n\le
x}\Ld(n)\lt(\lt[\f{x}{n}\rt]-2\lt[\f{x}{2n}\rt]\rt)+x\rt) \\
&=\cO\lt(\sum_{n\le
x}\Ld(n)\lt(\lt[\f{x}{n}\rt]-2\lt[\f{x}{2n}\rt]\rt)+x\rt) \\
&=\cO(x). \endsplit $$ From this, we have more generally the
following estimate
$$\pi\lt(\f{x}{2^k}\rt)\,\ln\f{x}{2^k}-\pi\lt(\f{x}{2^{k+1}}\rt)\,\ln\f{x}{2^{k+1}}=\cO\lt(\f{x}{2^k}\rt),\,\,k\in\BN.$$
Thus for any $K\in\BN$ we obtain that $$ \split \pi(x)\,\ln
x-\pi\lt(\f{x}{2^{K+1}}\rt)\,\ln\f{x}{2^{K+1}}&=\sum_{k=0}^K\lt(\pi\lt(\f{x}{2^k}\rt)\,\ln\f{x}{2^k}
-\pi\lt(\f{x}{2^{k+1}}\rt)\,\ln\f{x}{2^{k+1}}\rt) \\
&=\cO\lt(\sum_{k=0}^K\f{x}{2^k}\rt)=\cO(x). \endsplit $$ This
implies that $\,\ds\pi(x)=\cO\lt(\f{x}{\ln x}\rt).$ \qed

\proclaim{Proposition 5.5} The following asymptotic equation
$$\pi(x)\sim\f{x}{\ln x}$$ is equivalent to the asymptotic
equation $\,\psi(x)\sim x\,$ where the $\psi$-function is defined
by
$$\psi(x)=\sum_{n\le x}\Ld(n)=\sum_{p,\nu\ge 1 : p^{\nu}\le x}\ln
p\,.$$ $($ Here the function $\psi$ is introduced by Chebyshev.
$)$
\endproclaim

\noindent{\it Proof.} From the definition of the function $\psi$,
we have that $$\psi(x)=\sum_{p\le x}\lt[\f{\ln x}{\ln p}\rt]\,\ln
p\le\ln x\,\sum_{p\le x}\,1=\pi(x)\,\ln x.\tag{5.2}$$ On the other
hand, we note that for any $y$ with $1<y<x$,
$$ \split \pi(x)&=\pi(y)+\sum_{y<p\le x} 1\le\pi(y)+\sum_{y<p\le x}\f{\ln
p}{\ln y} \\ &\le c_2\cdot\f{y}{\ln y}+\f{\psi(x)}{\ln y}.
\endsplit$$ Thus, multiplying by the factor $\ln x/x$, the above
inequality becomes
$$\pi(x)\cdot\f{\ln x}{x}\le c_2\cdot\f{y\ln x}{x\ln
y}+\f{\psi(x)}{x}\cdot\f{\ln x}{\ln y}.\tag{5.3}$$ If we set
$y=x/\ln x$ in (5.3), then we have that
$$\pi(x)\cdot\f{\ln x}{x}\le\f{c_2}{\ln x-\ln\ln
x}+\f{\psi(x)}{x}\cdot\f{1}{1-\ds\f{\ln\ln x}{\ln
x}}\,.\tag{5.4}$$ Hence we complete the proof from (5.2) and
(5.4). \qed

\proclaim{Theorem 5.6[Mertens' Theorem]} If $p$ runs through all
prime numbers, then we have the following asymptotic
approximations;

$(a)$ $\ds\sum_{p\le x}\f{\ln p}{p}=\ln x+\cO(1),$ $\,\,(b)$
$\ds\sum_{p\le x}\f{1}{p}=\ln\ln x+c_3+\cO\lt(\f{1}{\ln x}\rt),$

$(c)$ $\ds\prod_{p\le x}\lt(1-\f{1}{p}\rt)=\f{c_4}{\ln
x}\lt(1+\cO\lt(\f{1}{\ln x}\rt)\rt),$

\noindent where $c_3>0$ and $c_4>0$ are some constants.
\endproclaim

\noindent{\it Proof.} (a) From (a) of Lemma 5.3 and Theorem 5.4,
we have that $$
\split x\ln x-x+\cO(\ln x)&=\sum_{n\le x}\Ld(n)\lt[\f{x}{n}\rt] \\
&=\sum_{p\le x}\lt[\f{x}{p}\rt]\,\ln p+\sum_{p\le\sqrt x,\nu\ge
2\,:\, p^{\nu}\le x}\lt[\f{x}{p^{\nu}}\rt]\,\ln p \\
&=\sum_{p\le x}\f{\ln p}{p}\cdot x-\sum_{p\le
x}\lt\{\f{x}{p}\rt\}\ln p+\cO\lt(\sum_{p\le\sqrt
x}\,\,\sum_{2\le\nu\le\f{\ln x}{\ln p}}\f{x}{p^{\nu}}\,\ln p\rt) \\
&=x\,\sum_{p\le x}\f{\ln p}{p}+\cO\lt(\sum_{p\le x}\ln
p\rt)+\cO\lt(x\,\sum_{n=1}^{\iy}\f{\ln n}{n^2}\rt) \\
&=x\,\sum_{p\le x}\f{\ln p}{p}+\cO\lt(\ln x\cdot c_2\cdot\f{x}{\ln x}\rt)+\cO\lt(x\,\sum_{n=1}^{\iy}\f{\ln n}{n^2}\rt) \\
&=x\,\sum_{p\le x}\f{\ln p}{p}+\cO(x). \endsplit $$ This implies
the first one.

(b) It follows from Proposition 4.2.2''[Abel transformation] that
$$ \split \sum_{p\le x}\f{1}{p}&=\sum_{p\le x}\f{\ln
p}{p}\cdot\f{1}{\ln p} \\ &=\f{1}{\ln x}\sum_{p\le x}\f{\ln
p}{p}+\int_2^x\sum_{p\le t}\f{\ln p}{p}\cdot\f{1}{t\ln^2 t}\,dt \\
&=1+\cO\lt(\f{1}{\ln x}\rt)+\int_2^x\f{1}{t\ln
t}\,dt+\int_2^x\lt(\sum_{p\le t}\f{\ln p}{p}-\ln t\rt)\f{1}{t\ln^2
t}\,dt.\endsplit $$ Since $\,a(t)=\ds\sum_{p\le t}\f{\ln p}{p}-\ln
t\,\,$ is bounded by (a), the following integral
$$\int_2^{\iy}\f{a(t)}{t\ln^2 t}\,dt$$ converges, and moreover we
have that $$\int_2^{\iy}\f{1}{t\ln t}\,dt=\ln\ln t-\ln\ln 2.$$
Therefore we conclude that $$ \split \sum_{p\le x}\f{1}{p}&=\ln\ln
x+\lt(1-\ln\ln 2+\int_2^{\iy}\f{a(t)}{t\ln^2
t}\,dt\rt)+\cO\lt(\f{1}{\ln x}+\int_x^{\iy}\f{|a(t)|}{t\ln^2
t}\,dt\rt) \\ &=\ln\ln x+c_3+\cO\lt(\f{1}{\ln x}\rt). \endsplit $$

(c) If we define the constant $c_5$ by
$$c_5=\sum_{n=2}^{\iy}\f{1}{n}\sum_p\f{1}{p^n},$$ then it follows
from simple calculation that $$ \split \ln\lt(\prod_{p\le
x}\lt(1-\f{1}{p}\rt)\rt)&=\sum_{p\le x}\ln
\lt(1-\f{1}{p}\rt)=-\sum_{p\le x}\sum_{n=1}^{\iy}\f{p^{-n}}{n} \\
&=-\sum_{p\le x}\f{1}{p}-\sum_{n=2}^{\iy}\f{1}{n}\sum_{p\le
x}\f{1}{p^n} \\
&=-\sum_{p\le
x}\f{1}{p}-c_5+\cO\lt(\,\sum_{n=2}^{\iy}\f{1}{n}\sum_{p>x}\f{1}{p^n}\rt)
\\ &=-\sum_{p\le
x}\f{1}{p}-c_5+\cO\lt(\,\sum_{n=2}^{\iy}\,\sum_{m>x}\f{1}{m^n}\rt) \\
&=-\sum_{p\le
x}\f{1}{p}-c_5+\cO\lt(\,\sum_{n=2}^{\iy}\f{1}{n}\cdot\f{1}{(n-1)
x^{n-1}}\rt) \\
&=-\sum_{p\le x}\f{1}{p}-c_5+\cO\lt(\f{1}{x}\rt). \endsplit $$
Hence this implies the required result. \qed

\proclaim{Lemma 5.7[Tauberian Theorem of Ingham and Newman]}

\noindent Let $F(t)$ be a bounded complex-valued function defined
on $(0,\iy)$ and integrable over every compact subset of
$(0,\iy)$, and let $G(z)$ be an analytic function defined on a
domain containing the closed half-plane $\Pi=\{z\in\BC :
\Rea(z)\ge 0\}$. If $G(z)$ agrees with the Laplace transformation
of $F(t)$ for all $z\in\Pi$, i.e. $$G(z)=\int_0^{\iy} F(t)\,e^{-z
t}\,dt,\,\,\Rea(z)>0,$$ then the improper integral $$\int_0^{\iy}
F(t)\,dt$$ converges. \endproclaim

\noindent{\it Proof.} Without loss of generality, we may assume
that $|F(t)|\le 1$ for all $t>0$. For $\ld>0$, we set
$$G_{\ld}(z)=\int_0^{\ld} F(t) e^{-zt}\,dt.$$ Then we see that
$G_{\ld}(z)$ is analytic on $\BC$. Thus it suffices to show that
$$\lim_{\ld\to\iy} G_{\ld}(0)=\lim_{\ld\to\iy}\int_0^{\ld}
F(t)\,dt=G(0).$$ Fix $\vep>0$. Then there are $\dt=\dt(\vep)>0$
and $R>0$ such that $1/R<\vep/3$ and $G(z)$ is analytic on the
compact region
$$\Om_{\dt,R}\fd\{z\in\BC : \Rea(z)\ge\dt, |z|\le R\}$$ with
boundary $\pa\Om_{\dt,R}=\gm$ which is a simple closed contour
oriented counterclockwise. By Cauchy integral formula, we have
that $$G(0)-G_{\ld}(0)=\f{1}{2\pi
i}\int_{\gm}\f{G(z)-G_{\ld}(z)}{z}\,dz.\tag{5.5}$$ We observe that
for $x=\Rea(z)>0$, $$|G(z)-G_{\ld}(z)|=\lt|\int_{\ld}^{\iy} F(t)
e^{-zt}\,dt\rt|\le\int_{\ld}^{\iy} e^{-xt}\,dt=\f{e^{-\ld
x}}{x},\tag{5.6}$$ and for $x=\Rea(z)<0$,
$$|G_{\ld}(z)|=\lt|\int_0^{\ld} F(t)
e^{-zt}\,dt\rt|\le\int_0^{\ld} e^{-xt}\,dt=\f{e^{-\ld
x}}{|x|}.\tag{5.7}$$ With technical reasons given in (5.6) and
(5.7), the relation (5.5) can be written again as
$$G(0)-G_{\ld}(0)=\f{1}{2\pi i}\int_{\gm}\,[G(z)-G_{\ld}(z)] e^{\ld
z}\,\lt(\f{1}{z}+\f{z}{R^2}\rt)\,dz.\tag{5.8}$$ If we denote by
$\gm_+$ the part of $\gm$ lying in $\Rea(z)>0$, then we see that
$$\f{1}{z} +\f{z}{R^2}=\f{2x}{R^2}$$ on $\gm_+$, and thus it
follows from (5.6) and (5.8) that $$ \split
|G(0)-G_{\ld}(0)|&\le\f{1}{2\pi}\int_{\gm_+}\lt|[G(z)-G_{\ld}(z)]
e^{\ld z}\lt(\f{1}{z}+\f{z}{R^2}\rt)\rt|\,dz \\
&\le\f{1}{2\pi}\cdot \f{e^{-\ld x}}{x}\cdot e^{\ld
x}\cdot\f{2x}{R^2}\cdot\pi R=\f{1}{R}<\f{\vep}{3}.\endsplit
\tag{5.9}$$ If we denote by $\gm_-$ the part of $\gm$ lying in
$\Rea(z)<0$, then we have that $$\f{1}{2\pi
i}\int_{\gm_-}\,G_{\ld}(z) e^{\ld
z}\,\lt(\f{1}{z}+\f{z}{R^2}\rt)\,dz=\f{1}{2\pi
i}\int_{|z|=R}\,G_{\ld}(z) e^{\ld
z}\,\lt(\f{1}{z}+\f{z}{R^2}\rt)\,dz$$ since $G_{\ld}(z)$ is
analytic on $\BC$. Thus similarly to (5.9) we obtain that
$$\lt|\f{1}{2\pi i}\int_{\gm_-}\,G_{\ld}(z) e^{\ld
z}\,\lt(\f{1}{z}+\f{z}{R^2}\rt)\,dz\rt|\le\f{1}{R}<\f{\vep}{3}.\tag{5.10}$$
Since the function $\ds G(z)\lt(\f{1}{z}+\f{z}{R^2}\rt)$ is
analytic on $\gm_-$, there is a constant $M=M(\dt,R)=M(\vep)>0$
such that $$\lt|G(z) e^{\ld z}\lt(\f{1}{z}+\f{z}{R^2}\rt)\rt|\le M
e^{\ld\Rea(z)}$$ for each $z\in\gm_-$. Since $\Rea(z)<0$ for
$z\in\gm_-$, the integral $$\f{1}{2\pi i}\int_{\gm_-}\,G(z) e^{\ld
z}\,\lt(\f{1}{z}+\f{z}{R^2}\rt)\,dz$$ tends to zero as
$\ld\to\iy$, and so there is a constant $N>0$ such that
$$\lt|\f{1}{2\pi i}\int_{\gm_-}\,G(z) e^{\ld
z}\,\lt(\f{1}{z}+\f{z}{R^2}\rt)\,dz\rt|<\f{\vep}{3}\tag{5.11}$$
whenever $\ld>N$. Thus if $\ld>N$, then it follows from (5.9),
(5.10), and (5.11) that $$|G(0)-G_{\ld}(0)|<\vep.$$ Therefore we
are done. \qed

\proclaim{Corollary 5.8[Simplified Version of the Theorem of
Weiner and Ikehara]}

Let $\,f(x)$ be a monotone nondecreasing function defined for
$x\ge 1$ with $f(x)=\cO(x)$. Suppose that $g(z)$ is analytic in
some region containing the closed half-plane $\Rea(z)\ge 1$ except
for a simple pole at $z=1$ with residue $\ap$ and, for any $z$
with $\Rea(z)>1$, $g(z)$ coincides with the Mellin transform of
$f(x)$, i.e. $$g(z)=z\int_1^{\iy}
f(x)\,x^{-z-1}\,dx,\,\,\Rea(z)>1.$$ Then we have that $f(x)\sim
\ap x$. \endproclaim

\noindent{\it Proof.} We note that the function $F(t)$ defined by
$$F(t)=e^{-t}\,f(e^t)-\ap$$ is bounded on $(0,\iy)$ and integrable
on each compact subset of $(0,\iy)$. Also its Laplace transform
$$ G(z)=\int_0^{\iy} [e^{-t}\,f(e^t)-\ap] e^{-zt}\,dt=\int_1^{\iy} f(x)\,x^{-z-2}\,dx-\f{\ap}{z}
=\f{1}{z+1}\,g(z+1)-\f{\ap}{z} \tag{5.12}$$ is well-defined in
$\Rea(z)>0$. By the assumption, the right-hand side of (5.12) is
analytic in some region containing the closed half-plane
$\Rea(z)\ge 0$. Thus it follows from Lemma 5.7 [Tauberian Theorem
of Ingham and Newman] that the improper integral
$$\int_0^{\iy}[e^{-t}\,f(e^t)-\ap]\,dt=\int_1^{\iy}\f{f(x)-\ap
x}{x^2}\,dx$$ converges. Now we shall prove that $f(x)\sim\ap x$
by using the nondecreasing monotonicity of $f$.

If $\,\ds\limsup_{x\to\iy}\f{f(x)}{x}>\ap$, then there exists some
$\dt>0$ so that $f(y)>(\ap+2\dt) y$ for infinitely many and
arbitrarily large $y$. Thus $f(x)>(\ap+2\dt) y>(\ap+\dt) x$ for
all $x$ with $y<x<\ds\lt(\f{\ap+2\dt}{\ap+\dt}\rt)\,y$, and
$$\int_y^{\lt(\f{\ap+2\dt}{\ap+\dt}\rt)\,y}\f{f(x)-\ap x}{x^2}\,dx>
\int_y^{\lt(\f{\ap+2\dt}{\ap+\dt}\rt)\,y}\f{\dt}{x}\,dx=\dt\cdot\ln\lt(\f{\ap+2\dt}{\ap+\dt}\rt)>0.$$
This gives a contradiction. So we conclude that
$$\limsup_{x\to\iy}\f{f(x)}{x}\le\ap.\tag{5.13}$$

If $\,\ds\liminf_{x\to\iy}\f{f(x)}{x}<\ap$, then there exists some
$\dt>0$ with $\dt<\ap/2$ so that $f(y)<(\ap-2\dt) y$ for
infinitely many and arbitrarily large $y$. Thus $f(x)<(\ap-2\dt)
y<(\ap-\dt) x$ for all $x$ with
$\ds\lt(\f{\ap-2\dt}{\ap-\dt}\rt)\,y<x<y$, and
$$\int^y_{\lt(\f{\ap-2\dt}{\ap-\dt}\rt)\,y}\f{f(x)-\ap x}{x^2}\,dx<
\int^y_{\lt(\f{\ap-2\dt}{\ap-\dt}\rt)\,y}\f{-\dt}{x}\,dx=-\dt\cdot\ln\lt(\f{\ap-\dt}{\ap-2\dt}\rt)<0.$$
This gives a contradiction. So we conclude that
$$\liminf_{x\to\iy}\f{f(x)}{x}\ge\ap.\tag{5.14}$$ Therefore we
complete the proof from (5.13) and (5.14). \qed

\proclaim{Lemma 5.9[Mertens]} $\zt(z)\neq 0$ for any $z$ with
$\Rea(z)=1$ and $z\neq 1$. \endproclaim

\noindent{\it Proof.} We observe that
$3+4\cos\th+\cos(2\th)=2(1+\cos\th)^2\ge 0$ for any $\th\in\BR$.
If $\zt(1+i t)=0$ for some $t\neq 0$, then the equation
$$\Th(s)={\zt(s)}^3\cdot {\zt(s+it)}^4\cdot\zt(s+2it)$$ has a zero
at $s=1$. Thus we have that $$\lim_{s\to 1}\ln
|\Th(s)|=-\iy.\tag{5.15}$$ Now it follows from Theorem 4.3.11 that
for any $s=\sm>1$, $$ \split
\ln |\zt(\sm+it)|&=-\Rea\lt(\sum_p\ln (1-p^{-\sm-it})\rt) \\
&=\Rea\lt(\sum_p\lt(p^{-\sm-it}+\f{1}{2}(p^2)^{-\sm-it}+\f{1}{3}(p^3)^{-\sm-it}+\cdots\rt)\rt)
\\ &=\Rea\lt(\,\sum_{n=1}^{\iy} b_n\,n^{-\sm-it}\rt) \endsplit $$
where $b_n$'s are certain nonnegative constants. This leads to the
following inequalities $$ \split \ln
|\Th(\sm)|&=\Rea\lt(\sum_{n=1}^{\iy}
b_n\,n^{-\sm}(3+4\,n^{-it}+n^{-2it})\rt) \\
&=\sum_{n=1}^{\iy} b_n\,n^{-\sm}(3+4\cos(t\,\ln n)+\cos(2t\,\ln
n))\ge 0,\endsplit $$ which contradict to (5.15). Hence we
complete the proof. \qed

\proclaim{Theorem 5.10[Prime Number Theorem]}

If $\,\pi(x)$ denotes the number of prime numbers $p\le x$, then
we have that $\ds\pi(x)\sim\f{x}{\ln x}$. \endproclaim

\noindent{\it Proof.} First of all, by Theorem 5.4[Chebyshev's
Theorem] we observe that $$ \split \psi(x)&=\sum_{p\le x}
\lt[\f{\ln x}{\ln p}\rt]\,\ln p\le\ln x\,\sum_{p\le x} 1 \\
&=\pi(x)\,\ln x=\cO(x).\endsplit $$ By Proposition 5.5, it
suffices to show that
$$\psi(x)\sim x.$$ By Theorem 4.3.18, the Mellin transform of $\psi(x)$ is
$$-\f{\zt'(z)}{\zt(z)}=z\int_1^{\iy}\f{\psi(x)}{x^{z+1}}\,dx,\,\,\Rea(z)>1.$$
In order to apply Corollary 5.8, we shall show that the function
$$-\f{\zt'(z)}{\zt(z)}-\f{1}{z-1}$$ is analytic in some region
containing the closed half-plane $\Rea(z)\ge 1$. By Proposition
4.3.16, there is some $\dt>0$ so that
$$\zt(z)=\f{1}{z-1}(1+h(z))$$ where $h(z)$ is analytic in
$B(1;\dt)$ and $|h(z)|<1$ there. Thus this implies that the
function $$-\f{\zt'(z)}{\zt(z)}-\f{1}{z-1}=-\f{h'(z)}{1+h(z)}$$ is
analytic at $z=1$. Finally, it follows from Proposition 4.3.16 and
Lemma 5.9 that the function $$-\f{\zt'(z)}{\zt(z)}-\f{1}{z-1}$$ is
analytic at any other points $z$ with $\Rea(z)=1$. Hence are are
done. \qed

\proclaim{Corollary 5.11} Let $f(x)$ be a number theoretic
function with nonnegative values and with $$\sum_{n\le x}
f(n)=\cO(x),$$ and let the Dirichlet series
$$F(z)=\sum_{n=1}^{\iy}\f{f(n)}{n^z}$$ be analytic in $\Rea(z)>1$
in the sense that the function $$F(z)-\f{\ap}{z-1}\,\,\text{ $($
$\ap$ is some fixed constant $)$}$$ is analytic in some region
containing the closed half-plane $\Rea(z)\ge 1$. Then we have that
$$\sum_{n\le x} f(n)\sim\ap x.$$ \endproclaim

\noindent{\it Proof.} It easily follows from Corollary 5.8 and the
following integral representation
$$F(z)=\sum_{n=1}^{\iy}\f{f(n)}{n^z}=z\int_1^{\iy}\lt(\sum_{n\le x}
f(n)\rt) x^{-z-1}\,dx.\,\,\,\,\,\,\,\,\,\,\,\qed$$

\proclaim{Corollary 5.12} Let $f(n)$ and $g(n)$ be two number
theoretic functions satisfying that $f(n)\ge 0$, $g(n)=\cO(f(n))$,
and $\ds\sum_{n\le x} f(n)=\cO(x)$. If two Dirichlet series
$$F(z)=\sum_{n=1}^{\iy}\f{f(n)}{n^z}\,\,\,\text{ and
}\,\,\,\,G(z)=\sum_{n=1}^{\iy}\f{g(n)}{n^z}$$ are analytic in
$\Rea(z)>1$ in the sense that the functions
$$F(z)-\f{\ap}{z-1},\,\,\,G(z)-\f{\bt}{z-1}\,\,\,\text{ $($ $\ap$ and
$\bt$ are some fixed constants $)$}$$ are analytic in some region
containing the closed half-plane $\Rea(z)\ge 1$, then we have that
$$\sum_{n\le x} g(n)\sim\gm x.$$ \endproclaim

\noindent{\it Proof.} First, we assume that $g(n)$ is real-valued.
Let us choose some constant $K>0$ so large that $|g(n)|\le K f(n)$
for all $n\in\BN$. We now apply Corollary 5.11 to the Dirichlet
series generated by the number theoretic function $h(n)=K
f(n)+g(n)$, given by $$H(z)=\sum_{n=1}^{\iy}\f{h(n)}{n^z}=K
F(z)+G(z).$$ By Corollary 5.11, we have that $$\sum_{n\le x}
h(n)=K\sum_{n\le x} f(n)+\sum_{n\le x} g(n)\sim K \ap x+\sum_{n\le
x} g(n)$$ and $$\sum_{n\le x} h(n)\sim K \ap x+\bt x.$$ This
implies the conclusion.

If $g(n)$ is complex-valued, then we set
$G^*(z)=\overline{G(\overline z)}$ and we consider
$$G_1(z)\fd\f{1}{2}[G(z)+G^*(z)]=\sum_{n=1}^{\iy}\f{\Rea(g(n))}{n^z}$$
and
$$G_2(z)\fd\f{1}{2i}[G(z)-G^*(z)]=\sum_{n=1}^{\iy}\f{\Ima(g(n))}{n^z}.$$
Hence we complete the proof by applying the above argument to
$G_1(z)$ and $G_2(z)$. \qed

In what follows, we furnish three examples as a foretaste of
importance of Corollary 5.12.

\proclaim{Corollary 5.13} If $\mu(n)$ is the M\"obius function and
$\ld(n)$ is the Liouville function, then we have that $$\sum_{n\le
x}\mu(n)=\fo(x)\,\,\,\text{ and }\,\,\,\sum_{n\le
x}\ld(n)=\fo(x).$$ \endproclaim

\noindent{\it Proof.} By Proposition 4.3.15, we apply Corollary
5.12 to the associated Dirichlet series $G(z)=1/\zt(z)$ and
$G(z)=\zt(2z)/\zt(z)$ which are analytic in some region containing
the closed half-plane $\Rea(z)\ge 1$. Since they have no
singularity at $z=1$, we conclude that $\bt=0$. \qed

As a third example, we consider the Dirichlet series
$$\zt_i(z)=\sum_{n=1}^{\iy}\f{r(n)}{n^z}$$ generated by the
number theoretic function $r(n)$ which counts the number of the
representations of $n$ as the sum of two squares. By Proposition
3.25 in Chapter 3, $r(n)$ can be considered as the number of
representations $n=\om\overline {\om}$ where $\om$ runs through
the ring $\BZ(i)$. Thus we obtain that
$$\zt_i(z)=\sum_{\om\in\BZ(i)\s\{0\}}\f{1}{|\om|^{2z}}
=\sum_{\om\in\BZ(i)\s\{0\}}\f{1}{(\om\overline\om\,)^z},$$ which
is called the {\bf $\zt$-function} for the number theory on the
ring $\BZ(i)$. In order to keep track of the arguments of
$\om\in\BZ(i)\s\{0\}$, Hecke originated the following Dirichlet
series
$$\Xi(h,z)=\sum_{\om\in\BZ(i)\s\{0\}}\f{1}{|\om|^{2z}}\cdot e^{4i
h\,\arg(\om)},\,\,h\in\BZ.$$ Then it is clear that
$\Xi(0,z)=\zt_i(z)$ and
$$\Xi(h,z)=\sum_{n=1}^{\iy}\f{1}{n^z}\lt(\sum_{|\om|^2 =n} e^{4i
h\,\arg(\om)}\rt), \,\,\,\Rea(z)>1.$$ Its convergence for
$\Rea(z)>1$ follows from the convergence of $\zt_i(z)$ for
$\Rea(z)>1$; which can be derived from the estimate $$\sum_{x\le
n\le
y}\f{r(n)}{n^z}=\f{1}{y^z}\,\cO(y-x)+z\int_x^y\cO(t-x)\,t^{-z-1}\,dt=\cO\lt(\f{1}{x^{z-1}}\rt)$$
which is obtained by applying Proposition 4.2.2[Abel
Transformation] and Proposition 4.2.8. The argument function
$\arg(\om)$ in $\Xi(h,z)$ is uniquely defined in
$-\pi<\arg(\om)\le\pi$.

\proclaim{Definition 5.14} Let $f$ be a complex-valued function
defined on $\BZ(i)$. Then $f$ is said to be multiplicative if
$f\not\equiv 0$ and $$f(\fm\fn)=f(\fm) f(\fn)\tag{5.16}$$ for any
pair $(\fm,\fn)\in\BZ(i)\times\BZ(i)$ with no common prime factor.
If $\,(5.16)$ holds for any pair $(\fm,\fn)\in\BZ(i)\times\BZ(i)$,
then we say that $f$ is completely multiplicative. \endproclaim

For instance, for $h\in\BZ$ we consider the function
$f(\om)=e^{4ih\,\arg(\om)}$. Then it is certainly completely
multiplicative and satisfies that $f(\fu)=1$ for unit elements
$\fu=1,i,-1,-i$. This is the reason why the factor $4$ in the
exponent was taken in $\Xi(h,z)$.

\proclaim{Proposition 5.15} Let $f$ be a complex-valued function
defined on $\BZ(i)$ satisfying that $f(\fu)=1$ for all units
$\fu\in\BZ(i)$. Suppose that the infinite series
$$F(z)=\sum_{\om\in\BZ(i)\s\{0\}}\f{f(\om)}{|\om|^{2z}}$$
converges absolutely for $\Rea(z)>\tau_0$.

$(a)$ If $f$ is multiplicative, then we have that for all $z$ with
$\Rea(z)>1$,
$$F(z)=4\prod_{\fp\in\BZ^+_p(i)}\lt(\,\sum_{\mu=1}^{\iy}\f{f(\fp^{\mu})}{|\fp|^{2\mu
z}}\rt)$$ where $\BZ^+_p(i)$ is the set of all prime elements
$\fp$ of $\BZ(i)$ with $0\le\arg(\fp)<\pi/2$.

$(b)$ If $f$ is completely multiplicative, then we have that for
all $z$ with $\Rea(z)>1$,
$$F(z)=4\prod_{\fp\in\BZ^+_p(i)}\f{1}{1-\ds\f{f(\fp)}{|\fp|^{2z}}}.$$

$(c)$ For $h\in\BZ$, we have that
$$\Xi(h,z)=4\prod_{\fp\in\BZ^+_p(i)}\f{1}{1-\ds\f{e^{4ih\,\arg(\fp)}}{|\fp|^{2z}}},\,\,\Rea(z)>1.$$
\endproclaim

\noindent{\it Proof.} It easily follows from the modification of
Proposition 4.3.13. \qed

\proclaim{Definition 5.16} We consider the function $\Ld_i$
defined on $\BZ(i)$ given by $$\Ld_i(\om)=\lt\{\aligned \ln
|\fp|,\,\,&\,\,\text{ if $\,\om=\fu\fp^{\nu}$
          for a unit $\fu$ and a prime $\fp$ } \\
                       0,  \,\,&\,\,\text{ if $\,\om$ is not such a
                       prime power, } \endaligned \rt.$$ which is
called the generalized Mangoldt function.  \endproclaim

In Chapter 4, we saw the relation between the Mangoldt function
and the quotient $\zt'(z)/\zt(z)$. Similarly, in what follows we
study the connection between the generalized Mangoldt function and
the quotient $$-\f{\Xi'(h,z)}{\Xi(h,z)};$$ in particular, this
quotient will play an important role in the Mellin transform of
the function
$$\psi_i(x)=\sum_{\om\in \BB_x(i)}\Ld_i(\om)\tag{5.17}$$ where
$\BB_x(i)=\{\om\in\BZ(i) : |\om|^2\le x\}$.

\proclaim{Lemma 5.17} For $\Rea(z)>1$ and $h\in\BZ$, we have that
$$-\f{\Xi'(h,z)}{\Xi(h,z)}=\f{1}{2}\sum_{\om\in\BZ(i)\s\{0\}}\f{\Ld_i(\om)}{|\om|^{2
z}}\,e^{4ih\,\arg(\om)}.$$ \endproclaim

\noindent{\it Proof.} Since $\log(1-e^{4ih\,\arg(\fp)}\cdot
|\fp|^{-2z})=\cO(|\fp|^{-2\Rea(z)}),$ the series $$H(z)\fd\log
4-\sum_{\fp\in\BZ^+_p(i)}\log\lt(1-\f{e^{4ih\,\arg(\fp)}}{|\fp|^{2z}}\rt)$$
converges uniformly in every compact subsets inside the half-plane
$\Rea(z)>1$, and so $H(z)$ is analytic in $\Rea(z)>1$. We also
have the relation $$e^{H(z)}=\Xi(h,z).$$ Thus we obtain that
$$H'(z)\cdot\Xi(h,z)=\Xi'(h,z).$$ Therefore we complete the proof
by calculating $H'(z)$ as follows;
$$ \split
H'(z)&=\sum_{\fp\in\BZ^+_p(i)}\f{1}{1-\ds\f{e^{4ih\,\arg(\fp)}}{|\fp|^{2z}}}\cdot
\f{e^{4ih\,\arg(\fp)}\cdot\log |\fp|^2}{|\fp|^{2z}} \\
&=2\sum_{\fp\in\BZ^+_p(i)}\f{\log |\fp|\cdot
e^{4ih\,\arg(\fp)}}{|\fp|^{2z}}\cdot\sum_{\mu=0}^{\iy}\f{e^{4ih\,\arg(\fp^{\mu})}}{|\fp^{\mu}|^{2z}}
\\ &=2\sum_{\fp\in\BZ^+_p(i)}\,\,\sum_{\mu=1}^{\iy}\f{\log |\fp|\cdot
e^{4ih\,\arg(\fp^{\mu})}}{|\fp^{\mu}|^{2z}} \\
&=\f{1}{2}\sum_{\fu\in\fU}\,\,\sum_{\fp\in\BZ^+_p(i)}\,\,\sum_{\mu=1}^{\iy}\f{\log
|\fu \fp|\cdot e^{4ih\,\arg((\fu \fp)^{\mu})}}{|(\fu
\fp)^{\mu}|^{2z}}
\\&=\f{1}{2}\sum_{\om\in\BZ(i)\s\{0\}}\f{\Ld_i(\om)}{|\om|^{2z}}\,e^{4ih\,\arg(\om)},
\endsplit $$ where $\fU$ denotes the set of all unit elements
$\fu$ of $\BZ(i)$. \qed

\proclaim{Lemma 5.18} For all $z$ with $\Rea(z)>1$, we have the
integral representation
$$-\f{\zt'_i(z)}{\zt_i(z)}=\f{z}{2}\int_1^{\iy}\f{\psi_i(x)}{x^{z+1}}\,dx$$
where $\psi_i$ is a function defined by
$\psi_i(x)=\sum_{\om\in\BB_x(i)}\Ld_i(\om).$ \endproclaim

\noindent{\it Proof.} By Lemma 5.17, we have that
$$-\f{\zt'_i(z)}{\zt_i(z)}=\f{1}{2}\sum_{\om\in\BZ(i)\s\{0\}}\f{\Ld_i(\om)}{|\om|^{2z}}.$$
It also follows from Proposition 4.2.2[Abel Transformation] that
$$\sum_{\om\in\BB_x(i)\s\{0\}}\f{\Ld_i(\om)}{|\om|^{2z}}=\f{1}{x^z}\cdot\psi_i(x)-\int_1^x\psi_i(y)\cdot
\f{-z}{y^{z+1}}\,dy.\tag{5.18}$$ From Proposition 3.24, we observe
that
$$\sum_{\fp\in\BZ_p(i),\,|\fp|^2\le
x}1\,\,\,\,\sim\,\,\,\,\pi(x)\tag{5.19}$$ where $\BZ_p(i)$ denotes
the set of all prime elements of $\BZ(i)$. Thus by the definition
of $\psi_i(x)$ and Theorem 5.4[Chebyshev's theorem] we obtain that
$$ \split
\psi_i(x)&=\sum_{\om\in\BB_x(i)}\Ld_i(\om)
=4\sum_{\fp\in\BZ_p(i),\,|\fp|^2\le x}\lt[\f{\ln x}{2\ln
|\om|}\rt]\,\ln |\om| \\
&=\cO\lt(\,\sum_{\fp\in\BZ_p(i),\,|\fp|^2\le x}\ln
x\rt)=\cO\lt(\ln x\sum_{\fp\in\BZ_p(i),\,|\fp|^2\le
x}1\rt)=\cO(x). \endsplit \tag{5.20}$$ Taking the limit $x\to\iy$
in (5.18), we can complete the proof. \qed

\proclaim{Lemma 5.19} For $h\in\BZ\s\{0\}$, we have that
$$\sum_{\om\in\BB_x(i)\s\{0\}} e^{4ih\,\arg(\om)}=\cO(\sqrt x\,\ln x). $$
\endproclaim

\noindent{\it Proof.} We write $\om=a+ib$ for $a,b\in\BZ$.
Observing that $\arg(a+ib)=\pi/2 -\arg(b+ia)$ for $a,b\in\BN$ and
considering only the sum over non-associated elements, we have
that $$\split \sum_{\om\in\BB_x(i)\s\{0\}}
e^{4ih\,\arg(\om)}&=4\sum_{a>0}\,\,\sum_{b\ge 0 : a^2 +b^2\le x}
e^{4ih\,\arg(a+ib)} \\
&=8\sum_{a>0}\,\,\sum_{b\ge a : a^2+b^2\le
x}\cos(4h\,\arg(a+ib))+\cO(\sqrt x) \\
&=8\sum_{0<a\le\sqrt {\f{x}{2}}}\,\,\sum_{a\le b\le\sqrt
{x-a^2}}\cos\lt(4h\tan^{-1}\lt(\f{b}{a}\rt)\rt)+\cO(\sqrt x).
\endsplit $$
Since $\tan^{-1}\lt(\f{\sqrt {x-a^2}}{a}\rt)-\tan^{-1} 1=\cO(1)$,
it follows from Proposition 4.2.3[The Euler Sum Formula] that
$$\split &\sum_{\om\in\BB_x(i)\s\{0\}} e^{4ih\,\arg(\om)} \\
&=8\sum_{1\le a\le\sqrt {\f{x}{2}}}\lt(\int_a^{\sqrt
{x-a^2}}\cos\lt(4h\tan^{-1}\lt(\f{y}{a}\rt)\rt)\,dy+\cO\lt(
1+\int_a^{\sqrt
{x-a^2}}\f{1}{a\lt(1+\f{y^2}{a^2}\rt)}\,dy\rt)\rt)+\cO(\sqrt x)
\\ &=8\sum_{1\le a\le\sqrt {\f{x}{2}}}\int_a^{\sqrt
{x-a^2}}\cos\lt(4h\tan^{-1}\lt(\f{y}{a}\rt)\rt)\,dy+\cO(\sqrt x) \\
&=8\int_1^{\sqrt {\f{x}{2}}}\int_t^{\sqrt {x-t^2}}\cos\lt(4h\tan^{-1}\lt(\f{y}{t}\rt)\rt)\,dy\,dt \\
&\qquad+\cO\lt(\sqrt x+\int_1^{\sqrt
{\f{x}{2}}}\lt|\f{d}{dt}\int_t^{\sqrt
{x-t^2}}\cos\lt(4h\tan^{-1}\lt(\f{y}{t}\rt)\rt)\,dy\rt|\,dt\rt)+\cO(\sqrt
x). \endsplit $$ We observe that $\ds\int_0^1\int_t^{\sqrt
{x-t^2}}\cos\lt(4h\tan^{-1}\lt(\f{y}{t}\rt)\rt)\,dy\,dt=\cO(\sqrt
x)\,$ and $$ \split \f{d}{dt}\int_t^{\sqrt
{x-t^2}}\cos\lt(4h\tan^{-1}\lt(\f{y}{t}\rt)\rt)\,dy&=\f{1}{t^2}\int_t^{\sqrt
{x-t^2}}\f{4h y\sin(\tan^{-1}(\f{y}{t}))}{1+\f{y^2}{t^2}}\,dy \\
&\quad -\f{t}{\sqrt {x-t^2}}\cos\lt(4h\tan^{-1}\lt(\f{\sqrt
{x-t^2}}{t}\rt)\rt)-\cos(h\pi) \\
&=\cO\lt( \int_t^{\sqrt {x-t^2}}\f{y}{t^2 +y^2}\,dy+\f{t}{\sqrt
{x-t^2}}\rt) \\
&=\cO\lt(\ln\lt(\f{x}{2 t^2}\rt)+1\rt).
\endsplit $$ Thus by applying polar coordinates $t=r\cos\th$ and $y=r\sin\th$
with $0<r\le\sqrt x$ and $\pi/4\le\th\le\pi/2$, we obtain that
$$ \split \sum_{\om\in\BB_x(i)\s\{0\}} e^{4ih\,\arg(\om)}&=8\int_0^{\sqrt
{\f{x}{2}}}\int_t^{\sqrt
{x-t^2}}\cos\lt(4h\tan^{-1}\lt(\f{y}{t}\rt)\rt)\,dy\,dt +\cO(\sqrt
x\,\ln x) \\ &=8\int_0^{\sqrt
x}\int_{\f{\pi}{4}}^{\f{\pi}{2}}\cos(4h\th)\,d\th \,r dr+\cO(\sqrt
x\,\ln x)=\cO(\sqrt x\,\ln x), \endsplit $$ because the last
integral vanishes for $h\in\BZ\s\{0\}$. Therefore we complete the
proof. \qed

\proclaim{Lemma 5.20} Let $f(n)$ be a number theoretic function
satisfying $$\lim_{N\to\iy}\f{1}{N}\sum_{n=1}^N f(n)=\ap.$$ For
$\Rea(z)>1$, we have the following formula
$$\sum_{n=1}^{\iy}\f{f(n)}{n^z}=\ap\cdot\zt(z)+\sum_{n=1}^{\iy}\lt(\f{1}{n^z}-\f{1}{(n+1)^z}\rt)
\lt(\,\sum_{m=1}^n f(m)-n\ap\rt).$$ \endproclaim

\noindent{\it Proof.} Applying Lemma 4.2.1[Abel Transformation],
we have that $$ \split &\sum_{n=1}^N\lt(\f{1}{(n+1)^z}
-\f{1}{n^z}\rt)\lt(\,\sum_{m=1}^n f(m)-\ap n\rt) \\
&\qquad=\sum_{n=1}^N\f{1}{(n+1)^z}\lt(\lt(\,\sum_{m=1}^{n+1}
f(m)-\ap(n+1)\rt)-\lt(\,\sum_{m=1}^n f(m)-\ap n\rt)\rt) \\
&\qquad\qquad-\f{1}{(N+1)^z}\lt(\,\sum_{m=1}^{N+1} f(m)-\ap
(N+1)\rt)+
(f(1)-\ap) \\
&\qquad=\sum_{n=1}^{N+1}\f{f(n)-\ap}{n^z}-\f{1}{(N+1)^z}\lt(\ap-\f{1}{N+1}\sum_{m=1}^{N+1}
f(m)\rt) \\
&\qquad=\sum_{n=1}^{N+1}\f{f(n)}{n^z}-\ap\sum_{n=1}^{N+1}\f{1}{n^z}+\cO\lt(\f{1}{(N+1)^{\Rea(z)-1}}\rt).
\endsplit $$ Since $(N+1)^{-(\Rea(z)-1)}$ tends to zero as
$N\to\iy$ for $\Rea(z)>1$, and also
$$ \split
\lt|\sum_{n=1}^N\lt(\f{1}{n^z}-\f{1}{(n+1)^z}\rt)\lt(\,\sum_{m=1}^n
f(m)-\ap n\rt)\rt|&=\lt|\sum_{n=1}^N z\lt(\int_n^{n+1}
\f{1}{x^{z+1}}\,dx\rt)\cdot n\lt(\ap-\f{1}{n}\sum_{m=1}^n
f(m)\rt)\rt|\\
&=\cO\lt(\,\sum_{n=1}^N\f{|z|}{n^{\Rea(z)}}\rt) \endsplit $$
converges for $\Rea(z)>1$, we can complete the proof by taking
$N\to\iy$. \qed

\proclaim{Lemma 5.21} For $h\in\BZ\s\{0\}$, $\Xi(h,z)$ has an
analytic continuation into the half-plane $\Rea(z)>1/2$.
Similarly, the function $$\zt_i(z)-\f{\pi}{z-1}$$ has an analytic
continuation into the half-plane $\Rea(z)>1/2$ in the sense that
$\zt_i(z)$ is analytic on $\Rea(z)>1/2$ except for a simple pole
at $z=1$ with residue $\pi$. \endproclaim

\noindent{\it Proof.} If we set $f(n)=\sum_{|\om|^2=n}
e^{4ih\,\arg(\om)}$ for $h\in\BZ\s\{0\}$, then it follows from
Lemma 5.19 that
$$\ap=\lim_{n\to\iy}\f{1}{N}\sum_{n=1}^N f(n)=0.$$ By Lemma 5.20,
we have that for $h\in\BZ\s\{0\}$,
$$\Xi(h,z)=\sum_{n=1}^{\iy}\lt(\f{1}{n^z}-\f{1}{(n+1)^z}\rt)\sum_{\om\in\BB_n(i)\s\{0\}}
e^{4ih\,\arg(\om)}.$$ Thus it follows from Lemma 5.19 that the
following sequence $$ \split
\sum_{n=M}^N\lt(\f{1}{n^z}-\f{1}{(n+1)^z}\rt)\sum_{\om\in\BB_n(i)\s\{0\}}
e^{4ih\,\arg(\om)}&=\cO\lt( |z|\sum_{n=M}^N\sqrt n\,\ln
n\lt|\int_n^{n+1}\f{1}{x^{z+1}}\,dx\rt|\rt) \\
&=\cO\lt( |z|\sum_{n=M}^N\f{\ln n}{n^{\Rea(z)+\f{1}{2}}}\rt)
\endsplit $$ converges uniformly to zero as $M\to\iy$ in every
compact subsets of the half-plane $\Rea(z)>1/2$. Hence this
implies the analytic continuation of $\Xi(h,z)$.

Similarly to the above, it follows from Lemma 5.20 that
$$ \split \zt_i(z)&=\sum_{n=1}^{\iy}\f{r(n)}{n^z}\\
&=\pi\cdot\zt(z)+\sum_{n=1}^{\iy}\lt(\f{1}{n^z}-\f{1}{(n+1)^z}\rt)\lt(\,\sum_{m=1}^n
r(m)-n\pi\rt). \endsplit $$ From Proposition 4.2.8, we see that
$$\sum_{m=1}^n r(m)-n\pi=\cO(\sqrt n).$$ Therefore we complete the
proof by applying the above argument once again. \qed

\proclaim{Lemma 5.22} For $h\in\BZ\s\{0\}$, $\Xi(h,z)\neq 0$ for
any $z$ with $\Rea(z)=1$. \endproclaim

\noindent{\it Proof.} It is trivial for the case $h=0$ and $z=1$,
because $\zt_i(z)$ has a pole at $z=1$. For the other cases, we
use a modified version of Lemma 5.9[Mertens].

Fix $h\in\BZ\s\{0\}$. If $\Xi(h,1+it)=0$ for some $t\neq 0$, then
the equation
$$\Th(s)={\zt_i(z)}^3\cdot {\Xi(h,s+it)}^4\cdot\Xi(2h,s+i2t)$$ has
a zero at $s=1$. Thus this implies that $$\lim_{s\to 1}\ln
|\Th(s)|=-\iy.\tag{5.21}$$ Now it follows from Proposition 5.15,
(c) that for any $s=\sm>1$,
$$ \split \ln |\Xi(h,\sm+it)|&=\ln 4-\sum_{\fp\in\BZ_p^+(i)}\ln
\lt|1-\f{e^{4ih\,\arg(\fp)}}{|\fp|^{2\sm+i2t}}\rt| \\
&=\ln 4+\sum_{\fp\in\BZ_p^+(i)}\,\sum_{n=1}^{\iy}\f{\cos n
(4h\,\arg(\fp)-2t\ln |\fp|)}{n |\fp|^{2n\sm}}. \endsplit $$ This
leads to the following inequalities $$\ln |\Th(s)|=8\ln
4+\sum_{\fp\in\BZ_p^+(i)}\,\sum_{n=1}^{\iy}\f{3+4\cos n
(4h\,\arg(\fp)-2t\ln |\fp|)+\cos n(8h\,\arg(\fp)-4t\ln |\fp|)}{n
|\fp|^{2n\sm}}\ge 0,$$ which contradicts to (5.21). Hence we
complete the proof. \qed

\proclaim{Proposition 5.23}
$\,\,\ds\psi_i(x)\fd\sum_{\om\in\BB_x(i)}\Ld_i(\om)\,\,\sim\,\,
2x.$\endproclaim

\noindent{\it Proof.} It is trivial that $\psi_i(x)$ is a monotone
non-decreasing function on $[0,\iy)$. By (5.20), we have
$\,\psi_i(x)=\cO(x)$. Thus it follows from Lemma 5.18 and Lemma
5.21 that the function $-\zt'_i(z)/\zt_i(z)$ given by
$$-\f{\zt'_i(z)}{\zt_i(z)}=z\int_1^{\iy}\f{1}{2}\psi_i(x)\,\f{1}{x^{z+1}}\,dx$$
is analytic in $\Rea(z)>1$ and the function
$$-\f{\zt'_i(z)}{\zt_i(z)}-\f{1}{z-1}$$ has an analytic
continuation into some region containing the closed half-plane
$\Rea(z)\ge 1$. Therefore Corollary 5.11 implies the conclusion.
\qed

\proclaim{Proposition 5.24} $\,\,\ds\sum_{\om\in\BB_x(i)}
e^{4ih\,\arg(\om)}\,\Ld_i(\om)=\fo(x)\,\,$ for $h\in\BZ\s\{0\}$.
\endproclaim

\noindent{\it Proof.} We observe that
$e^{4ih\,\arg(\om)}\,\Ld_i(\om)=\cO(\Ld_i(\om))$ for
$h\in\BZ\s\{0\}$ and $\om\in\BZ(i)\s\{0\}$. From Lemma 5.17 and
Lemma 5.22, two Dirichlet series
$$-\f{\zt'(z)}{\zt(z)}=\f{1}{2}\sum_{\om\in\BZ(i)\s\{0\}}\f{\Ld_i(\om)}{|\om|^{2
z}}\,\,\,\text{ and }\,\,\,
-\f{\Xi'(h,z)}{\Xi(h,z)}=\f{1}{2}\sum_{\om\in\BZ(i)\s\{0\}}\f{\Ld_i(\om)\,e^{4ih\,\arg(\om)}}{|\om|^{2
z}},\,\,\,h\in\BZ\s\{0\},$$ are analytic in $\Rea(z)>1$ and have
an analytic continuation with no singularity at $z=1$ into some
region containing the closed half-plane $\Rea(z)\ge 1$. Therefore
Corollary 5.12 and Proposition 5.23 imply the required one. \qed

\proclaim{Theorem 5.25[Hecke's Prime Number Theorem for the ring
$\BZ(i)$]}

$(a)$ If $\,\pi_i(x)$ denotes the number of all non-associated
prime elements $\fp$ with $|\fp|^2\le x$, i.e. the number of all
prime elements in $\BZ_p^+(i)\cap\BB_x(i)$, then we have that
$$\pi_i(x)\,\,\sim\,\,\f{x}{\ln x}.$$

$(b)$ If $\,\pi_i(x;\ap,\bt)$ denotes the number of all prime
elements $\fp\in\BZ_p(i)\cap\BB_x(i)$ with $\ap\le\arg(\fp)<\bt$
for $0\le\ap<\bt\le 2\pi$, then we have that
$$\pi_i(x;\ap,\bt)\,\,\sim\,\,\f{2}{\pi}\,(\bt-\ap)\,\f{x}{\ln
x}.$$ \endproclaim

\noindent{\it Proof.} We observe the following estimate
$$\split \sum_{k\ge 2}\,\,\sum_{\fp\in\BZ_p^+(i) :|\fp|^{2k}\le x}
e^{4ih\,\arg(\fp^k)}\,\ln|\fp|&=\cO\lt(\ln x\sum_{k\ge
2}\,\,\sum_{\fp\in\BZ_p^+(i) :|\fp|^2\le x^{1/k}}\,1\rt) \\
&=\cO\lt(\ln x \sum_{2\le k\le\ln x/\ln 2}\f{\sqrt x}{\ln\sqrt
x}\rt) \\ &=\cO(\sqrt x\,\ln x). \endsplit $$ This implies that
$$\split 4\sum_{\fp\in\BZ_p^+(i)\cap\BB_x(i)}
e^{4ih\,\arg(\fp)}\,\ln|\fp|&=4\sum_{k\ge
1}\,\,\sum_{\fp\in\BZ_p^+(i) :|\fp|^{2k}\le x}
e^{4ih\,\arg(\fp^k)}\,\ln|\fp|+\cO(\sqrt x\,\ln x)\\
&=\sum_{\om\in\BB_x(i)} e^{4ih\,\arg(\om)}\,\Ld_i(\om)+\cO(\sqrt
x\,\ln x) \\
&=\lt\{\aligned 2 x+\fo(x), \,\,\,& h=0, \\
          \fo(x),\,\,\, & h\neq 0. \endaligned \rt. \endsplit $$
Thus it follows from the above estimate and Proposition 4.2.2
[Abel Transformation] that
$$\split \sum_{\fp\in\BZ_p^+(i)\cap\BB_x(i)}
e^{4ih\,\arg(\fp)}&=\sum_{\fp\in\BZ_p^+(i) : 2\le |\fp|^2\le x}
e^{4ih\,\arg(\fp)}\,\ln |\fp|^2\cdot\f{1}{\ln |\fp|^2} \\
&=\f{1}{\ln x}\sum_{\fp\in\BZ_p^+(i) : 2\le |\fp|^2\le x}
e^{4ih\,\arg(\fp)}\,\ln |\fp|^2 \\
&\qquad-\int_2^x\sum_{\fp\in\BZ_p^+(i) : 2\le |\fp|^2\le t}
e^{4ih\,\arg(\fp)}\,\ln
|\fp|^2\cdot\f{-1}{t\ln^2 t}\,dt\\
&=\f{2}{\ln x}\sum_{\fp\in\BZ_p^+(i)\cap\BB_x(i)}
e^{4ih\,\arg(\fp)}\,\ln|\fp| +\cO\lt(\int_2^x\f{1}{\ln^2
t}\,dt\rt) \\
&=\lt\{\aligned \f{x}{\ln x}+\fo\lt(\f{x}{\ln
x}\rt)+\cO\lt(\f{x}{\ln x}\rt),\,\,\, & h=0, \\
\fo\lt(\f{x}{\ln x}\rt)+\cO\lt(\f{x}{\ln x}\rt),\,\,\, & h\neq 0
\endaligned \rt. \\
&=\lt\{\aligned \f{x}{\ln x}+\fo\lt(\f{x}{\ln
x}\rt),\,\,\, & h=0, \\
\fo\lt(\f{x}{\ln x}\rt),\,\,\, & h\neq 0. \endaligned \rt.
\endsplit \tag{5.22}$$

(a) By (5.22) on $h=0$, we have that
$$\pi_i(x)=\sum_{\fp\in\BZ_p^+(i)\cap\BB_x(i)}
e^{4ih\,\arg(\fp)}=\f{x}{\ln x}+\fo\lt(\f{x}{\ln x}\rt).$$

(b) It easily follows from (5.22) on $h\neq 0$ that
$$\split \lim_{x\to\iy}\f{1}{4\pi_i(x)}\sum_{\fp\in\BZ_p(i)\cap\BB_x(i)}
e^{2\pi
ih(\f{2}{\pi}\,\arg(\fp))}&=\lim_{x\to\iy}\f{1}{4\pi_i(x)}\cdot 4
\sum_{\fp\in\BZ_p^+(i)\cap\BB_x(i)} e^{2\pi ih(\f{2}{\pi}\,\arg(\fp))} \\
&=\lim_{x\to\iy}\f{1}{\pi_i(x)}\sum_{\fp\in\BZ_p^+(i)\cap\BB_x(i)}
e^{2\pi ih(\f{2}{\pi}\,\arg(\fp))} \\
&=\lim_{x\to\iy}\f{\ln x}{x}\sum_{\fp\in\BZ_p^+(i)\cap\BB_x(i)}
e^{4ih\,\arg(\fp)}=0.\endsplit $$ Thus by Theorem 2.13 [Weyl's
Criterion] we see that the sequence
$$\{\th_{\fp,x}\fd\f{2}{\pi}\,\,\arg(\fp) :
\fp\in\BZ_p(i)\cap\BB_x(i),\,x\in\BR_+\}$$ is uniformly
distributed modulo $2\pi$. Hence by Proposition 2.12 we have that
$$ \lim_{x\to\iy}\f{1}{4\pi_i(x)}\sum_{\fp\in\BZ_p(i)\cap\BB_x(i)}
f(\arg(\fp))=\f{1}{2\pi}\int_0^{2\pi} f(\th)\,d\th\tag{5.23}$$ for
any real-valued Riemann integrable function $f(\th)$ on
$[0,2\pi)$. If we take $f(\th)=\chi_{[\ap,\bt)}(\th)$ in (5.23),
we obtain that $$\lim_{x\to\iy}\f{\pi_i(x;\ap,\bt)}{4\pi_i(x)}=
\lim_{x\to\iy}\f{1}{4\pi_i(x)}\sum_{\fp\in\BZ_p(i)\cap\BB_x(i)}
f(\arg(\fp))=\f{1}{2\pi}\int_0^{2\pi}
f(\th)\,d\th=\f{1}{2\pi}(\bt-\ap).$$ Therefore this implies the
required result. \qed

\enddocument